\documentclass[12pt,a4paper,twoside]{article}

\newcommand{\pageformat}[6]{\setlength{\hoffset}{-1in}
                  \setlength{\voffset}{-1in}
                  \addtolength{\hoffset}{#5}
                            \addtolength{\voffset}{#6}
                            \setlength{\oddsidemargin}{#1}
                            \setlength{\evensidemargin}{#2}
                            \setlength{\textwidth}{\paperwidth}
                  \addtolength{\textwidth}{-\oddsidemargin}
                  \addtolength{\textwidth}{-\evensidemargin}
                  \addtolength{\textwidth}{-\marginparsep}
                  \addtolength{\textwidth}{-\marginparwidth}
                            \setlength{\topmargin}{#3}
                            \setlength{\textheight}{\paperheight}
                  \addtolength{\textheight}{-\topmargin}
                  \addtolength{\textheight}{-\headheight}
                  \addtolength{\textheight}{-\headsep}
                  \addtolength{\textheight}{-\footskip}
                  \addtolength{\textheight}{-#4}}
\pageformat{2cm}{3cm}{25mm}{25mm}{1pt}{0pt}

\usepackage{ifthen}
\newboolean{article}
    \setboolean{article}{true}
\newboolean{report}
\newboolean{book}
\newboolean{letter}
\newboolean{german}
\newboolean{italian}
\newboolean{nobaselinestretch}
\newboolean{nosectionappendix}
\newboolean{oldtoc}
\newboolean{nosectionequn}
\newboolean{notheorem}

\ifthenelse{\boolean{german}}{
    \usepackage{german}}{}

\usepackage[latin1]{inputenc}

\usepackage{amsmath}
\usepackage{amssymb}
\usepackage[mathscr]{eucal}

\ifthenelse{\boolean{notheorem}}{}{
    \usepackage{theorem}}

\ifthenelse{\boolean{nobaselinestretch}}{}{
    \renewcommand{\baselinestretch}{1.25}}

\newenvironment{env}[2]{\begin{#1}#2\end{#1}}{}
    \newcommand{\beq}[1]{\begin{env}{equation}{#1}}
    \newcommand{\beqn}[1]{\begin{env}{equation*}{#1}}
    \newcommand{\bal}[1]{\begin{env}{align}{#1}}
    \newcommand{\baln}[1]{\begin{env}{align*}{#1}}
    \newcommand{\bga}[1]{\begin{env}{gather}{#1}}
    \newcommand{\bgan}[1]{\begin{env}{gather*}{#1}}
    \newcommand{\bflal}[1]{\begin{env}{flalign}{#1}}
    \newcommand{\bflaln}[1]{\begin{env}{flalign*}{#1}}
    \newcommand{\bmu}[1]{\begin{env}{multline}{#1}}
    \newcommand{\bmun}[1]{\begin{env}{multline*}{#1}}
    \newcommand{\bsp}[1]{\begin{env}{split}{#1}}

    \newcommand{\eeq}{\end{env}}
    \newcommand{\eeqn}{\end{env}}
    \newcommand{\eal}{\end{env}}
    \newcommand{\ealn}{\end{env}}
    \newcommand{\ega}{\end{env}}
    \newcommand{\egan}{\end{env}}
    \newcommand{\eflal}{\end{env}}
    \newcommand{\eflaln}{\end{env}}
    \newcommand{\emu}{\end{env}}
    \newcommand{\emun}{\end{env}}
    \newcommand{\esp}{\end{env}}

\newcommand{\lf}{\vspace{2ex}}

\renewcommand{\bf}[1]{\textbf{#1}}
\renewcommand{\it}[1]{\textit{#1}}

\renewcommand{\sf}[1]{\textsf{#1}}

\renewcommand{\tt}[1]{\texttt{#1}}
\newcommand{\hl}[1]{\bf{\it{#1}}}

\newcommand{\msf}[1]{\text{\small $\sf{#1}$}}

\newcommand{\cmc}[1]{\mathcal{#1}}
\newcommand{\eus}[1]{\mathscr{#1}}

\newcommand{\bb}[1]{\mathbb{#1}}

\newcommand{\nbd}[1]{$#1$\nobreakdash--}
\newcommand{\ol}[1]{\overline{#1}}

\newcommand{\wt}[1]{\widetilde{#1}}

\newcommand{\bnorm}[1]{\bigl\lVert#1\bigr\rVert}

\newcommand{\Bnorm}[1]{\Bigl\lVert#1\Bigr\rVert}

\newcommand{\bfam}[1]{\bigl(#1\bigr)}
\newcommand{\Bfam}[1]{\Bigl(#1\Bigr)}
\newcommand{\AB}[1]{\langle#1\rangle}

\newcommand{\CB}[1]{\{#1\}}
\newcommand{\bCB}[1]{\bigl\{#1\bigr\}}
\newcommand{\BCB}[1]{\Bigl\{#1\Bigr\}}
\newcommand{\SB}[1]{[#1]}
\newcommand{\bSB}[1]{\bigl[#1\bigr]}

\newcommand{\LO}[1]{(#1]}

\newcommand{\RO}[1]{[#1)}

\newcommand{\set}[2][]{
    \ifthenelse{\equal{#1}{}}{
        \CB{#2}}{
        \CB{#1~|~#2}}}
\newcommand{\bset}[2][]{
    \ifthenelse{\equal{#1}{}}{
        \bCB{#2}}{
        \bCB{#1~|~#2}}}
\newcommand{\Bset}[2][]{
    \ifthenelse{\equal{#1}{}}{
        \BCB{#2}}{
        \BCB{#1~\big|~#2}}}
\newcommand{\zero}{\CB{0}}

\DeclareMathOperator{\ls}{\normalfont\msf{span}}

\DeclareMathOperator{\id}{\normalfont\msf{id}}

\newcommand{\C}{\bb{C}}

\newcommand{\N}{\bb{N}}

\newcommand{\R}{\bb{R}}

\newcommand{\cB}{\cmc{B}}
\newcommand{\cC}{\cmc{C}}

\newcommand{\cI}{\cmc{I}}

\newcommand{\sB}{\eus{B}}

\newcommand{\sS}{\eus{S}}
\newcommand{\sT}{\eus{T}}

\newcommand{\I}{{I\!\!\!\;I}}

\newcommand{\s}{\text{\scriptsize$\sS$}}
\newcommand{\st}{\text{\scriptsize$\sT$}}

\ifthenelse{\boolean{nosectionequn}}{}{
    \numberwithin{equation}{section}
    }

\ifthenelse{\boolean{article}\or\boolean{letter}\or\boolean{nosectionequn}}{
    \setboolean{nosectionappendix}{true}}{}
\ifthenelse{\boolean{nosectionappendix}}{}{
    \renewcommand{\appendix}{
        \chapter*{\appendixname}
        \addcontentsline{toc}{chapter}{\appendixname}
        \renewcommand{\thesection}{\Alph{section}}
        \setcounter{section}{0}}}

\ifthenelse{\boolean{report}\or\boolean{book}}{
    }{}

\ifthenelse{\boolean{notheorem}}{}{
        \newcommand{\mnname}{Mathematical note.}
        \newcommand{\enname}{End of the note.}
        \newcommand{\definame}{Definition.}
        \newcommand{\propname}{Proposition.}
        \newcommand{\lemname}{Lemma.}
        \newcommand{\exname}{Example.}
        \newcommand{\exername}{Exercise.}
        \newcommand{\remname}{Remark.}
        \newcommand{\obname}{Observation.}
        \newcommand{\thmname}{Theorem.}
        \newcommand{\corname}{Corollary.}
        \newcommand{\proofname}{Proof.}
    \ifthenelse{\boolean{german}}{
        \renewcommand{\mnname}{Mathematische Notiz.}
        \renewcommand{\enname}{Ende der Notiz.}
        \renewcommand{\exname}{Beispiel.}
        \renewcommand{\exername}{Übung.}
        \renewcommand{\remname}{Bemerkung.}
        \renewcommand{\obname}{Beobachtung.}
        \renewcommand{\thmname}{Satz.}
        \renewcommand{\corname}{Korollar.}
        \renewcommand{\proofname}{Beweis.}}{}
    \ifthenelse{\boolean{italian}}{
        \renewcommand{\mnname}{Nota matematica.}
        \renewcommand{\enname}{Fina della nota.}
        \renewcommand{\definame}{Definizione.}
        \renewcommand{\propname}{Proposizione.}
        \renewcommand{\exname}{Esempio.}
        \renewcommand{\exername}{Esercizio.}
        \renewcommand{\remname}{Nota.}
        \renewcommand{\obname}{Osservazione.}
        \renewcommand{\thmname}{Teorema.}
        \renewcommand{\corname}{Corollario.}
        \renewcommand{\proofname}{Dimostrazione.}

       \renewcommand{\appendixname}{Appendice}

       }{}
    \theoremheaderfont{\normalfont\bfseries}
    \theoremstyle{change}
        \theorembodyfont{\rmfamily}
            \newtheorem{emp}{}[section]
                \newcommand{\bemp}[1][]{
                    \begin{emp}\hskip-\labelsep\bf{#1}\hskip\labelsep}
                \newcommand{\eemp}{\end{emp}}
\newtheorem{itemp}[emp]{}
                \newcommand{\bitemp}[1][]{
                    \begin{itemp}\hskip-\labelsep\bf{#1}\hskip\labelsep\normalfont\itshape}
                \newcommand{\eitemp}{\end{itemp}}
            \newtheorem{mn}[emp]{\mnname}
                \newcommand{\bnm}{\begin{mn}~\begin{quotation}\renewcommand{\baselinestretch}{1}\small\noindent\ignorespaces}
                \newcommand{\enm}{\end{quotation}\hfill\bf{\enname}\end{mn}}
            \newtheorem{ex}[emp]{\exname}
                \newcommand{\bex}{\begin{ex}}
                \newcommand{\eex}{\end{ex}}
            \newtheorem{exer}[emp]{\exername}
                \newcommand{\bexer}{\begin{exer}}
                \newcommand{\eexer}{\end{exer}}
            \newtheorem{defi}[emp]{\definame}
                \newcommand{\bdefi}{\begin{defi}}
                \newcommand{\edefi}{\end{defi}}
            \newtheorem{rem}[emp]{\remname}
                \newcommand{\brem}{\begin{rem}}
                \newcommand{\erem}{\end{rem}}
            \newtheorem{ob}[emp]{\obname}
                \newcommand{\bob}{\begin{ob}}
                \newcommand{\eob}{\end{ob}}
        \theorembodyfont{\normalfont\itshape}
            \newtheorem{thm}[emp]{\thmname}
                \newcommand{\bthm}{\begin{thm}}
                \newcommand{\ethm}{\end{thm}}
            \newtheorem{prop}[emp]{\propname}
                \newcommand{\bprop}{\begin{prop}}
                \newcommand{\eprop}{\end{prop}}
            \newtheorem{cor}[emp]{\corname}
                \newcommand{\bcor}{\begin{cor}}
                \newcommand{\ecor}{\end{cor}}
            \newtheorem{lem}[emp]{\lemname}
                \newcommand{\blem}{\begin{lem}}
                \newcommand{\elem}{\end{lem}}
\newenvironment{empn}[1]{\lf\noindent\bf{#1}\ignorespaces\hskip\labelsep}{\lf}
        \newcommand{\bempn}[1]{\begin{empn}{#1}}
        \newcommand{\eempn}{\end{empn}}
        \newcommand{\bitempn}[1]{\begin{empn}{#1}\normalfont\itshape}
        \newcommand{\eitempn}{\end{empn}}
                \newcommand{\bnmn}{\begin{empn}{\mnname}~\begin{quotation}\renewcommand{\baselinestretch}{1}\small\noindent\ignorespaces}
                \newcommand{\enmn}{\end{quotation}\hfill\bf{\enname}\end{empn}}
        \newcommand{\bexn}{\begin{empn}{\exname}}
        \newcommand{\eexn}{\end{empn}}
        \newcommand{\bexern}{\begin{empn}{\exername}}
        \newcommand{\eexern}{\end{empn}}
        \newcommand{\bdefin}{\begin{empn}{\definame}}
        \newcommand{\edefin}{\end{empn}}
        \newcommand{\bremn}{\begin{empn}{\remname}}
        \newcommand{\eremn}{\end{empn}}
        \newcommand{\bobn}{\begin{empn}{\obname}}
        \newcommand{\eobn}{\end{empn}}

\newcommand{\qedsymbol}{~\rule[-0.35mm]{2mm}{2mm}}
    \newcounter{proof}[emp]
    \newenvironment{Proof}[1]{
        \vspace{1ex}
        \renewcommand{\item}[1][\stepcounter{proof}(\roman{proof})]%
            {##1\hskip\labelsep}
        \noindent\textsc{#1\hskip\labelsep}}{
        \nolinebreak\qedsymbol}
    \newcommand{\proof}[1][\proofname]{
        \begin{Proof}{#1}\ignorespaces}
    \newcommand{\qed}{\end{Proof}}
    \newcommand{\noqed}{
        \renewcommand{\qedsymbol}{}
        \end{Proof}}}
    \ifthenelse{\boolean{italian}}{
        \renewcommand{\proofname}{Dimostrazione.}}{}

\usepackage[varg]{txfonts}

\usepackage{hyperref}

\begin{document}

\bibliographystyle{amsalpha}

\title{Pure semigroups of isometries on Hilbert $C^*$-modules}

\author{B.V.Rajarama Bhat\thanks{Statistics and Mathematics Unit, Indian Statistical Institute Bangalore, R.\ V.\ College Post, Bangalore 560059, India. E-mail: \href{mailto:bhat@isibang.ac.in}{\tt{bhat@isibang.ac.in}}}  ~and Michael Skeide\thanks{Dipartimento E.G.S.I., Università degli Studi del Molise, Via de Sanctis, 86100 Campobasso, Italy. E-mail \href{mailto:skeide@unimol.it}{\tt{skeide@unimol.it}}}}

\date{August 2014}

\maketitle

\begin{abstract}
\noindent 
We show that pure strongly continuous semigroups of adjointable isometries on a Hilbert $C^*$-module are standard right shifts. By counter examples, we illustrate that the analogy of this result with the classical result on Hilbert spaces by Cooper, cannot be improved further to understand arbitrary semigroups of isometries in the classical way. The counter examples include a strongly continuous semigroup of non-adjointable isometries, an extension of the standard right shift that is not strongly continuous, and a strongly continuous semigroup of adjointable isometries that does not admit a decomposition into a maximal unitary part and a pure part.
\end{abstract}

\section{Introduction}

The following classical result characterizes one-parameter semigroups (in the sequel, semigroups) of isometries on Hilbert spaces; see the standard text book Sz.-Nagy and Foias \cite[Theorem 9.3]{Sz-NFo70} or Sz.-Nagy \cite{Sz-Na64} for a short self-contained treatment.

\bitemp[Theorem (Cooper \cite{Coo47}).] \label{SzNFthm}
Let $\s=\{ \s_t\}_{t\in\R+}$ be a strongly continuous semigroup
of isometries on a Hilbert space $H$. Then $H=H_u\oplus H_p$ where $H_u$ and $H_p$ are unique subspaces such that:
\begin{enumerate}
\item
$H_u$ reduces $\s$ to a semigroup of unitaries.

\item
$H_p$ reduces $\s$ to a completely nonunitary semigroup.
\end{enumerate}
Moreover, the completely nonunitary part is unitarily equivalent to the standard right shift on $L^2(\R_+,K)$ for some multiplicity Hilbert space $K$.
\eitemp

It is our aim to prove the appropriately formulated analogue of the last sentence for semigroups of isometries on Hilbert modules (Theorem \ref{mainthm}, proved in Section \ref{proofSEC}) and to analyze to what extent we can save statements about decomposition (Section \ref{(c)eSEC}).

\noindent
Recall that:
\begin{itemize}
\item
A family $\s=\{ \s_t\}_{t\in\R+}$ of linear maps on a vector space is a (one-parameter) \hl{semigroup} if $\s_r\s_t=\s_{r+t}$ and if $\s_0=\id$, the identity.

\item
A semigroup $\s$ on a normed space $V$ is \hl{strongly continuous} if $t\mapsto\s_tv$ is continuous for all $v\in V$.

\item
A linear map $\s_t$ on a (pre-)Hilbert module $E$ (for instance, on a (pre-)Hilbert space) is an \hl{isometry} if it preserves inner products: $\AB{\s_tx,\s_tx}=\AB{x,y}$ for all $x,y\in E$. (It follows that a semigroup of isometries is strongly continuous if and only if $t\mapsto\AB{x,\s_ty}$ is continuous for all $x,y\in E$.) An isometry between Hilbert space is, like all bounded linear operators, \hl{adjointable}, that is, it has an adjoint. An isometry between (pre-)Hilbert modules is adjointable if and only if there exists a projection onto its range, that is, if and only if its range is \hl{complemented}.

\item
A \hl{unitary} is a surjective isometry. A semigroup $\s$ of isometries is \hl{completely nonunitary} if there is no nonzero subspace that reduces $\s$ to a semigroup of unitaries; see below.  A semigroup $\s$ of adjointable isometries is \hl{pure} if $\s_t^*$ converges strongly to $0$ for $t\to\infty$. (Note that, by the Riemann-Lebesgue lemma, the unitary semigroup of multiplication operators $e^{i\lambda t}$ on $L^2(\R)$ converges to $0$ weakly.) Equivalently, the projections $\s_t\s_t^*$ converge strongly to $0$, respectively, the projections $\id-\s_t\s_t^*$ converge strongly to $\id$. For a semigroup $\s$ of isometries on Hilbert space, the property to be completely nonunitary and the property to be pure are equivalent. See, however, Example \ref{nondecex}.

\item
A subspace $W$ of a vector space $V$ is \hl{invariant} for the semigroup $\s$ if $\s_tW\subset W$ for all $t$. We say $W$ \hl{reduces} $\s$ and the (co)restriction of $\s$ to $W$ is called the semigroup \hl{reduced} by $W$ or simply the \hl{reduced semigroup}. (Note that a semigroup of isometries is completely nonunitary if and only if $\s_tW=W$ for all $t$ implies $W=\zero$.) A (pre-)Hilbert submodule $F$ of a (pre)Hilbert module $E$ is \hl{completely reducing}%
\footnote{%
This property is known in literature as \hl{reducing subspace}. In order to avoid confusion with invariant subspaces that reduce a semigroup, we prefer an explicitly different terminology.
}
for the semigroup $\s$ on $E$ if both $F$ and $F^\perp$ (see Section \ref{(c)eSEC}) reduce $\s$. There are two special instances of completely reducing submodules: Firstly, if $\s$ is adjointable and $F$ reduces also the adjoint semigroup $\s^*$. (In general, if $\s$ reduces $F$, then $\AB{\s^*_t(F^\perp),F}=\AB{F^\perp,s_tF}\subset\AB{F^\perp,s_tF}=\zero$, so $\s_t^*(F^\perp)\subset F^\perp$, that is, $F^\perp$ reduces $\s^*$. So if $F$ reduces $\s^*$, then $F^\perp$ reduces $\s=\s^{**}$.) Secondly, if $F$ reduces the semigroup of isometries $\s$ to a unitary semigroup. (Indeed, from $\s_tF=F$ and $\s_t$ being an isometry, we obtain $\AB{\s_t(F^\perp),F}=\AB{\s_t(F^\perp),\s_tF}=\AB{F^\perp,F}=\zero$, so $\s_t(F^\perp)\subset F^\perp$.)

For Hilbert spaces there are many equivalent characterizations of their closed completely reducing subspaces, and each completely reducing closed subspace decomposes the semigroup into a direct sum. (Already for pre-Hilbert spaces that latter statement is false; and Hilbert modules behave much more like pre-Hilbert spaces than like Hilbert spaces.) But only the definition we give here, meets the situations we deal with in Section \ref{(c)eSEC}, when we analyze how much of the first part of Theorem \ref{SzNFthm} can be saved.

\item
Let $F$ be a Hilbert module over a $C^*$-algebra $\cB$. The \hl{standard right shift} over $F$, called \hl{multiplicity module}, is the semigroup $v=\CB{v_t}_{t\in\R_+}$ of isometries on $L^2(\R_+,F):=L^2(\R_+)\otimes F$ defined by
\beqn{
\SB{v_tf}(x)
~=~
\begin{cases}
f(x-t)&x\ge t,
\\
0&\text{else.}
\end{cases}
}\eeqn
Here, the \hl{external tensor product} $L^2(\R_+)\otimes F$ of the Hilbert \nbd{\C}module $L^2(\R_+)$ and the Hilbert \nbd{\cB}modules $F$ can be identified with the completion of the space of functions $\ls\bCB{x\mapsto h(x)y\,|\,h\in\cC_c(\R_+),y\in F}$ ($C_c$ meaning continuous functions with compact support) in the norm arising from the inner product $\AB{f,g}:=\int\AB{f(x),g(x)}\,dx$. Clearly, the the standard right shift is strongly continuous. (Indeed, for any bounded strongly continuous  functions $t\mapsto a_t\in\sB(L^2(\R_+))$ the function $t\mapsto a_t\otimes\id_F\in\sB^a(L^2(\R_+)\otimes F)$ is strongly continuous, too, because by boundedness it is sufficient to check strong continuity on the total set of elementary tensors, and because on elementary tensors strong continuity of $a_t\otimes\id_F$ follows from strong continuity of $a_t$.)
\end{itemize}

\noindent
In Section \ref{proofSEC}, we will prove:
\begin{thm}\label{mainthm}
Let $\{ \s_t\}_{t\geq 0}$ be a pure strongly continuous semigroup
of adjointable isometries on a Hilbert module $E$.  Then $\{
\s_t\}$ is unitarily equivalent to the standard right shift on $L^2(\R_+,F)$ over some multiplicity module $F$.
\end{thm}

The reader who is interested only in that result may switch to the proof in Section \ref{proofSEC}, immediately. In the remainder of this introduction, we motivate the proof. We also explain briefly why for Hilbert modules we need a fresh proof. In the Section \ref{(c)eSEC} we explain by counter examples, why Theorem \ref{mainthm}, as compared with Theorem \ref{SzNFthm}, is the best we may hope for. We also explain why we may not hope for a generalization of the Stone-von Neumann theorem to Hilbert modules, and why for von Neumann (or $W^*$-)modules no new proof is needed.

\paragraph{Motivation of the proof.}
Suppose we knew that $E=L^2(R_+,F)$ ``in some way'' (that is, after choosing a suitable unitary $V\colon L^2(R_+,F)\rightarrow E$) and that, in this identification, $\s$ is the standard right shift (that is, $\s_t=Vv_tV^*$ under the unitary equivalence transform arising from that suitable unitary).

The first problem to be faced is, how can we extract $F$ from the abstract space $E$ and the abstract isometry semigroup $\s$? From the solution we propose here, all the other steps will suggest themselves.

Note that $E$ contains many copies of $F$, namely, for each $0\le a<b$ the submodule $\I_{\RO{a,b}}F$ of functions $\I_{\RO{a,b}}y$ (where $y\in F$ and $\I_S$ denotes the indicator function of the set $S$), is isomorphic to $F$ via $\I_{\RO{a,b}}y\mapsto y\sqrt{b-a}$. Moreover, the isomorphism
\beqn{
L^2(R_+,F)
~\supset~
\I_{\RO{a,b}}F
~\longrightarrow~
F
~\longrightarrow~
\I_{\RO{a+t,b+t}}F
~\subset~
L^2(R_+,F)
}\eeqn
we obtained that way, is nothing but $\s_t$ restricted to $\I_{\RO{a,b}}F$. On the other hand, the elements $\I_{\RO{a,b}}y$ (varying also $a$ and $b$) form a total subset of $E$. So, understanding how to get them abstractly will, if Theorem \ref{mainthm} is true, lead to a proof.

First of all, assuming the isometries are adjointable, we easily get the submodules $E_{a,b}:=L^2(\RO{a,b},F)$ of $F$. Indeed, the projection $\s_t\s_t^*$ onto the range of $\s_t$, is nothing but multiplication with $\I_{\RO{t,\infty}}$, and the projection onto $E_{a,b}$ is
\beqn{
\I_{\RO{a,b}}
~=~
\I_{\RO{a,\infty}}-\I_{\RO{b,\infty}}
~=~
\s_a\s_a^*-\s_b\s_b^*
~=:~
p_{a,b}.
}\eeqn
So, how to find, inside $E_{a,b}$, the elements of the form $\I_{\RO{a,b}}y$? We define $u^{a,b}_t$ to be the unitary shift modulo $b-a$ on $E_{a,b}$. (This can be done entirely in terms of $\s$; see \eqref{uabdefi}.) Then $\I_{\RO{a,b}}y$ are precisely those elements that are invariant under that unitary group. We find them by applying the projection
\beqn{
q_{a,b}
\colon
x
~\longmapsto~
\frac{1}{b-a}\int_a^bu^{a,b}_tx\,dt
}\eeqn
to all elements $x\in E_{a,b}$. (Thanks to strong continuity of $\s$, this integral is a well-defined Riemann integral over a continuous vector-valued function.) In this way, we get elements in $E_{a,b}$ that behave like elements $\I_{\RO{a,b}}y$. The analytic heart of the proof will be to show that the elements $q_{c,d}x$ ($x\in E_{a,b}$, $\RO{c,d}\subset\RO{a,b}$) are total in $E_{a,b}$. This is done in Lemma \ref{limit}, which asserts that
\beqn{
\sum_{k=1}^nq_{\frac{k-1}{n},\frac{k}{n}}
}\eeqn
converges strongly to the identity of $E_{0,1}$. The rest is simple reconstruction of the right shift on $E=\bigoplus_{k\in\N}E_{k-1,k}\cong E_{0,1}^\infty$ out of $u^{0,1}_t$. (See Lemma \ref{interlem}.)

\lf
The key point in the proof is that we manage (basically by Lemma \ref{limit}) to approximate explicitly in norm an arbitrary element of $E$ by sums over elements that behave like $\I_{\RO{a,b}}y$. No arguments like zero-complement or weakly total=strongly total (which work only for Hilbert spaces) are involved. The original proof for Hilbert spaces involves unbounded operators and adjoints; it does not appear to be generalizable to modules.

\newpage

\section{(Counter)examples and other obstacles}\label{(c)eSEC}

The proofs of several statements in the classical Theorem \ref{SzNFthm} for Hilbert spaces rely on several crucial properties and results that are not available for Hilbert modules. The most important are: Self-duality of closed subspaces and, therefore, existence of projections onto them; existence of adjoints; weakly total subsets are norm total. In this section, we explain the consequences of having these pieces missing by counter examples and prove some weaker statements. The proof of that statement we can confirm fully, Theorem \ref{mainthm}, has been outlined in the end of the introduction and will be performed in Section \ref{proofSEC}.

Closely related to projections onto closed submodules are the subtleties around orthogonal complements. For convenience, we start by repeating those facts that generalize easily to Hilbert modules. (Of course, orthogonal complements can be defined for arbitrary subsets, and a number of statements remain true also for pre-Hilbert modules and their not necessarily closed submodules, for instance, such as the \nbd{\cB}linear span of a subset. We ignore these, here.)

\begin{itemize}
\item
Let $F$ and $G$ be closed submodules of a Hilbert \nbd{\cB}modules $E$ such that $\AB{F,G}=\zero$. Then $y+z\mapsto y\oplus z$ defines an isometry $F+G\rightarrow F\oplus G$. This isometry is, clearly, surjective, that is, a unitary. We denote this situation as $F+G=F\oplus G$. Of course, $F\cap G=\zero$.

\item
Let $F$ be a closed submodule of a Hilbert \nbd{\cB}module $E$. The \hl{orthogonal complement} of $F$ is defined as $F^\perp:=\CB{x\in E\colon\AB{F,x}=\zero}$. Since, clearly, $\AB{F,F^\perp}=\zero$, we have $F+F^\perp= F\oplus F^\perp$. (Corollary: $F\oplus F^\perp=E$ if and only if the exists a projection $p\in\sB^a(E)$ with $pE=F$.) Obviously, $F_1\subset F_2$ implies $F_1^\perp\supset F_2^\perp$.

\item
Clearly, $F^{\perp\perp}\supset F$. Applying this to $F^\perp$ we get $(F^\perp)^{\perp\perp}\supset F^\perp$, and by the preceding conclusion we get $(F^{\perp\perp})^\perp\subset F^\perp$. So, $F^{\perp\perp\perp}=F^\perp$. (Corollary: If $G=F^\perp$ for some submodule $F$, then $G^{\perp\perp}=G$. Corollary: An adjointable map $a\colon E\rightarrow E'$ is zero on $F^{\perp\perp}$ if and only if it is zero on $F$. Indeed, $aF=\zero$ $\Leftrightarrow$ $\AB{E',aF}=\zero$ $\Leftrightarrow$ $\AB{a^*E',F}=\zero$ $\Leftrightarrow$ $a^*E'\subset F^\perp=(F^{\perp\perp})^\perp$ $\Leftrightarrow$ $\AB{a^*E',F^{\perp\perp}}=\zero$ $\Leftrightarrow$ $\AB{E',aF^{\perp\perp}}=\zero$ $\Leftrightarrow$ $aF^{\perp\perp}=\zero$. Note, too, that it appears to be unknown if this statement is true for all bounded right linear maps $a$. If this was true, one could show that if $G$ is a submodule of $E$ containing $F$ and fulfilling $F^\perp\cap G=\zero$, then $G^\perp=F^\perp$ and, consequently, $G^{\perp\perp}=F^{\perp\perp}\supset G$.)

\item
Since also $F^\perp+F^{\perp\perp}=F^\perp\oplus F^{\perp\perp}$, we get $(F\oplus F^\perp)^\perp=F^\perp\cap F^{\perp\perp}=\zero$. So, $F$ and $F^\perp$ \hl{separate} the points of $E$ and the \hl{relative orthogonal complement} of $F$ in $F^{\perp\perp}$ is $\zero$, that is, $F$ separates the points of $F^{\perp\perp}$. Clearly, $F^{\perp\perp}$ is the biggest submodule containing $F$ and orthogonal to $F^\perp$ (that is, subset of $F^{\perp\perp}$). So, if we have orthogonal submodules $F$ and $G$ of $E$ and if $G=F^\perp$, then $F^{\perp\perp}\oplus F^\perp$ is the biggest submodule of $E$ that allows the decomposition into a direct sum with summands containing $F$ and $G$, respectively. (Once more, if the values of bounded right linear maps on $F^{\perp\perp}$ should turn out to be determined uniquely by the values on $F$, it would be possible to show that for two orthogonal submodules $F$ and $G$ separating the points of $E$ ($F^\perp\cap G^\perp=\zero$), $F^{\perp\perp}\oplus G^{\perp\perp}$ is the unique maximal choice to embed $F\oplus G$ into a direct sum contained in $E$. Without that, it is only easy to see that also $F^{\perp\perp}$ and $G^{\perp\perp}$ are orthogonal.)

\end{itemize}

\noindent
As the most general situation, we consider a semigroup $\s$ of isometries on a Hilbert module
 $E$ that is \it{a priori} neither adjointable nor strongly continuous. Example \ref{nonadex} illustrates, among many peculiarities regarding decomposition, that a strongly continuous semigroup of isometries need not be adjointable. In Observation \ref{notscob} we point out that a `surprisingly reasonable' semigroup of adjointable isometries (namely, an extension of a standard right shift) is not strongly continuous.

\bex \label{nonadex}
Let $\cB$ be a $C^*$-algebra and let $\cI$ be a proper closed ideal of $\cB$. We mention the following technical result: For any Hilbert space $H$, the complement of $H\otimes\cI$ in $H\otimes\cB$ is $(H\otimes\cI)^\perp=H\otimes(\cI^\perp)$ (Indeed, by means of an ONB $(e_s)_{s\in S}$ for $H$, the elements of $H\otimes\cB$ are precisely those of the form $x=\sum_se_e\otimes b_s$ where the sum $\sum_sb_s^*b_s$ exists. Of course, if all $b_s$ are in $\cI^\perp$, then $\AB{H\otimes\cI,x}=\zero$. If one $b_s$ is not $\cI^\perp$, then there exists $c\in\cI$ such that $c^*b_s\ne0$. Consequently, $\AB{e_s\otimes c,x}=c^*b_s\ne0$, so that $x$ is not in the complement of $H\otimes\cI$.) Consequently, $H\otimes\cI$ has zero-complement in $H\otimes\cB$ if and only if $\cI$ is essential in $\cB$, and $H\otimes\cI$ is complemented in $H\otimes\cB$ if and only if $\cI$ is complemented in $\cB$.

Consider the the bilateral right shift $u_t$ on $L^2(\R,\cB)$, obviously, a strongly continuous semigroup. Clearly, $u_t$ sends $E:=L^2(\R_-,\cI)\oplus L^2(\R_+,\cB)$ into $E$, so that the the (co)restrictions $\s_t$ of $u_t$ to $E$ define a strongly continuous semigroup of isometries. An isometry is adjointable if and only if its image is complemented. So, if $\cI$ is not complemented in $\cB$, then $\s_tE=L^2(\R_-,\cI)\oplus L^2(\SB{0,t},\cI)\oplus L^2(\RO{t,\infty},\cB)$ is not complemented in $E=L^2(\R_-,\cI)\oplus L^2(\SB{0,t},\cB)\oplus L^2(\RO{t,\infty},\cB)$, so that $\s_t$ is not adjointable.

If $\cI$ is essential in $\cB$, then even $(\s_tE)^\perp=\zero$ for all $t$. It would be an interesting question to examine all strongly continuous semigroups of (not necessarily) adjointable isometries with this property. (If $(\s_tE)^\perp=\zero$ and $\s_t$ is adjointable, then $\s_t$ is a unitary.)

On the other hand, if $\cI$ is complemented in $\cB$ (so that $\s_t$ is adjointable), then $E=L^2(\R,\cI)\oplus L^2(\R_+,\cI^\perp)$. On the first summand, $\s_t$ reduces to a (unitary) bilateral  right shift. On the second summand it reduces to a standard right shift over $\cI^\perp$.
\eex

What about the decomposition stated in Theorem \ref{SzNFthm} in the general situation? Well, let us first speak about the maximal unitary part. Note that $\s_tE$ is a decreasing family of (closed) submodules of $E$. Clearly, the (closed!) submodule $E_u:=\bigcap_{t\in\R_+}\s_tE$ reduces the semigroup $\s_t$. Since $\s_tE$ is decreasing, the restriction is a unitary onto $E_u$. Moreover, if $E'$ is any other closed submodule of $E$ (completely) reducing $\s_t$ to a semigroups of unitaries, then $E'=\s_tE'\subset\s_tE$ so that $E'\subset E_u$. In other words, $E_u$ carries the unique maximal unitary part of $\s$.

In Example \ref{nonadex}, we get $E_u=L^2(\R,\cI)$, and the maximal unitary part is complementary if and only if $\cI$ is. Moreover, $E_u^\perp=L^2(\R_+,\cI^\perp)$ reduces $\s$ to a pure semigroup of adjointable isometries (already given as a standard right shift). We have $E=E_u\oplus E_u^\perp$ if and only if $\cI$ is complemented in $\cB$. In particular, if $\cI$ is essential (and proper), then $E_u^\perp$ is $\zero$ but $E_u\ne E=E_u^{\perp\perp}$.

What about $E_u^\perp$ in the general situation? Well, since $E_u$ reduces $\s$ to a unitary semigroup, $E_u$ is even completely reducing, that is, also $E_u^\perp$ reduces $\s$. We must have $\bigcap_{t\in\R_+}\s_t(E_u^\perp)=\zero$. (Otherwise, this subset of $\s_0(E_u^\perp)=E_u^\perp$ would contribute to $\bigcap_{t\in\R_+}\s_tE=E_u$.) So, $E_u^\perp$ reduces $\s$ to a completely nonunitary semigroup.

Note that, in Example \ref{nonadex}, the submodule $E':=L^2(\R_+,\cI)$ of $E_u$ fulfills $\bigcap_{t\in\R_+}\s_tE'=\zero$. So, also $E_u^\perp\oplus E'=L^2(\R_+,\cI\oplus\cI^\perp)$ reduces $\s$ to a completely nonunitary semigroup. There is no such condition as  ``the maximal submodule that reduces $\s_t$ to a completely nonunitary semigroup'' that can replace pureness. However, $(E_u^\perp\oplus E')^\perp=L^2(\R_-,\cI)$, which does not reduce $\s$. So, $E_u^\perp\oplus E'$ is not completely reducing. We do not know if $E_u^\perp$ is something like the biggest completely reducing submodule module that reduces $\s$ to a completely nonunitary semigroup. (The main problem is to show that if $F_1$ and $F_2$ reduce $\s$ to completely nonunitary semigroups, then so does $\ol{F_1+F_2}$.)

We collect:

\bprop\label{uperop}
A semigroup $\s$ of isometries on $E$ (neither necessarily adjointable nor necessarily strongly continuous) has a unique maximal unitary part acting on the completely reducing submodule $E_u$ and $E_u^\perp$ reduces $\s$ to a completely nonunitary semigroup. $E_u\oplus E_u^\perp$ need not be all of $E$.
\eprop

Note that it is unclear whether $E_u^{\perp\perp}$ is always invariant, so that both $E_u^\perp$ and $E_u^{\perp\perp}$ would be completely reducing, or not. (In Example \ref{nonadex}, $E_u^{\perp\perp}=L^2(\R_-,\cI)\oplus L^2(\R_+,\cI^{\perp\perp})$ is invariant for all choices of $\cI$. But in an attempt to prove that this is true for all semigroups of isometries, we again bump into the question if a bounded right linear map is determined by its values on a submodule with zero-complement.) However, if $E_u^{\perp\perp}\ne E_u$ is reduces $\s$, then by maximality of $E_u$, the semigroup reduced to $E_u^{\perp\perp}\ne E_u$ is certainly not unitary.

There is an interesting submodule of $E_u^\perp$, namely, $E_p:=\ol{\bigcup_t(\s_tE)^\perp}$. (Indeed, since $\s_tE$ is not smaller than $E_u$, the complement $(\s_tE)^\perp$ is not bigger than the complement $E_u^\perp$, and this turns over to union and closure.) This submodule is interesting, because if $\s$ is pure so that $\id_E-\s_t\s_t^*$ converges strongly to $\id_E$, then $E_p=\ol{\bigcup_t(\s_tE)^\perp}=\ol{\bigcup_t(\id_E-\s_t\s_t^*)E}=E$. And whenever the statements of Theorem \ref{SzNFthm} hold in full, then $E_u^\perp=E_p$.

Generally, also the submodule $E_p$ reduces $\s$ to a completely nonunitary semigroup of isometries. (Indeed, one easily verifies that $\s_t((\s_r E)^\perp)\subset (\s_{t+r}E)^\perp$, so that $\s_tE_p=\ol{\bigcup_r(\s_t(\s_rE)^\perp})\subset E_p$. Of course, a submodule $E_u^\perp$ can reduce only to a completely nonunitary semigroup.) In Example \ref{nonadex}, we have $E_p=E_u^\perp$, but in Example \ref{nondecex} we will see that this need not be so, not even even if $\s$ is adjointable. In general, we also do not know, if $E_p^\perp=E_u$. (Note, however, that if $E_p^\perp\supset E_u$ reduces $\s$ to a semigroup of unitaries, then necessarily $E_p^\perp=E_u$, because $E_u$ is maximal.) This situation improves if $\s$ is adjointable.

In Example \ref{nonadex}, the restriction of $\s$ to $E_u^\perp=E_p$ is a a standard right shift. In particular, the restricted $\s_t$ are adjointable independently on whether the original $\s_t$ were adjointable or not. Also here, we do not know if this is true in general for one of the restrictions of $\s$ to $E_u^\perp$ or to $E_p$. However, if the $\s_t$ are adjointable from the beginning, then
\beqn{
E_p^\perp
~=~
\Bfam{\bigcup_t(\s_tE)^\perp}^\perp
~=~
\bigcap_t(\s_tE)^{\perp\perp}
~=~
\bigcap_t\s_tE,
}\eeqn
because $\s_tE$ is complemented, so, $E_p^\perp=E_u$ and, further, $E_u^{\perp\perp}=E_p^{\perp\perp\perp}=E_p^\perp=E_u$. Also, the restriction of $\s_t$ to $E_p$ remains adjointable. (Suppose $x\in(\s_rE)^\perp$, so that $\AB{x,\s_ry}=0$ for all $y\in E$. So, $\AB{\s_t^*x,\s_ry}=\AB{x,\s_r(\s_ty)}=0$, that is, $\s_t^*x\in(\s_rE)^\perp$. It follows that the adjoint $\s_t^*$ leaves $E_p$ invariant, and the restriction of $\s_t^*$ is an adjoint of the restricted $\s_t$.) Finally, the restriction to $E_p$ is pure. (Indeed, to show that $\s_t^*y\to0$ for all $y\in E_p$, by boundedness of $\s_t^*$ it is sufficient to show that for $x$ from the dense subset $\bigcup_r(\s_rE)^\perp$. So let $y\in(\s_rE)^\perp$ for some $r$. Then $0=\AB{\s_rx,y}=\AB{x,\s_r^*y}$ for all $x\in E$, hence, $\s_r^*y=0$. It follows that $\s_ty=0$ for all $t\ge r$.) Moreover, $E_p$ is the unique biggest completely reducing submodule of $E_u^\perp$ that reduces $s$ to a pure semigroup of isometries.

We collect:

\bprop\label{pprop}
A semigroup $\s$ of adjointable isometries on $E$ (not necessarily strongly continuous) has a unique maximal pure part acting on $E_p$ such that $E_p$ is completely reducing and such that $E_p^\perp=E_u$ (so that $E_p^\perp$ reduces $\s$ to the unique maximal unitary part). $E_u\oplus E_p$ need not be all of $E$ (that is, $E_p^{\perp\perp}=E_u^\perp$ need not be $E_p$) but (since $E_u=E_p^\perp$) its orthogonal complement is $\zero$.

Clearly, if $\s$ is also strongly continuous, then to the part on $E_p$ we may apply Theorem \ref{mainthm}.
\eprop

We now discuss an example that shows that even assuming adjointability, the module $E$ need not coincide with $E_u\oplus E_p$. Since in this example $E_u=\zero$, it also follows that completely nonunitary ($E_u=\zero$) does not imply pure ($E_p=E$). We prepare with the following lemma (the second part of which will also be important in the proof in Section \ref{proofSEC}).

\blem\label{interlem}
Let $u_t$ denote the unitary right-shift modulo $1$ on $L^2\RO{0,1}$. Suppose $\breve{\s}$ is an isometry on a Hilbert \nbd{\cB}module $\breve{E}$. For $t\in\R_+$ denote by $n_t$ the largest integer $\le t$.
\begin{enumerate}
\item\label{1}
The maps
\beqn{
\s_t
~:=~
(u_t\otimes\id_{\breve{F}})
\bfam{\I_{\RO{0,1-(t-n_t)}}\otimes\breve{\s}^{n_t}+\I_{\RO{1-(t-n_t),1}}\otimes\breve{\s}^{n_t+1}}
}\eeqn
define a strongly continuous semigroup of adjointable isometries on $L^2\RO{0,1}\otimes\breve{E}$.

\item\label{2}
If $\breve{E}=F^\infty$ and if $\breve{\s}\colon(y_1,y_2,\ldots)\mapsto(0,y_1,y_2,\ldots)$ is the one-sided right shift, then $\s_t$, under the canonical isomorphism
\beqn{
L^2\RO{0,1}\otimes F^\infty
~\cong~
L^2\RO{0,1}^\infty\otimes F
~\cong~
L^2(\R_+)\otimes F,
}\eeqn
is nothing but the standard right shift $v_t$ on $L^2(\R_+)\otimes F$.
\end{enumerate}
\elem

\proof
Clearly, $\s_t$ are adjointable isometries. The function $u_t$ is strongly continuous and the projections $\I_{\RO{0,1-(t-n_t)}}$ and $\I_{\RO{1-(t-n_t),1}}$ in $\sB(L^2\RO{0,1})$ depend strongly continuously on $t\in\RO{n-1,n}$ for all $n\in\N$. So, for the same reason for which the standard right shift $v_t$ is strongly continuous as discussed before Theorem \ref{mainthm}, also $\s_t$ is strongly right continuous. On the other hand, the definition of $\s_t$ does not change (for $t>0$), if we replace $n_t$ with the largest integer $<t$. So, $\s_t$ is also strongly left continuous. For showing Part \ref{1} it, therefore, remains to show that the $\s_t$ form a semigroup.

Note that $\s_n=\id\otimes\breve{\s}^n$ for all $n\in\N_0$. (In particular, $\s_0=\id$.) Also, $\s_t\s_n=\s_{t+n}=\s_n\s_t$ for all $t\in\R_+$ and $n\in\N_0$, is easily verified. So, it is sufficient to verify the semigroup property $\s_r\s_t=\s_{r+t}$ for $r,t\in(0,1)$. Note that
\baln{
\I_{\RO{1-r,1}}u_t
~=~
&
u_t
\begin{cases}
\I_{\RO{1-r-t,1-t}}								&t\le1-r,
\\
\I_{\RO{0,1-t}}+\I_{\RO{2-r-t,1}}		&t\ge1-r,
\end{cases}
\intertext{hence,}
\I_{\RO{0,1-r}}u_t
~=~
(\id-\I_{\RO{1-r,1}})u_t
~=~
&
u_t
\begin{cases}
\I_{\RO{0,1-r-t}}+\I_{\RO{1-t,1}}		&t\le1-r,
\\
\I_{\RO{1-t,2-r-t}}								&t\ge1-r.
\end{cases}
}\ealn
Therefore,
\baln{
\s_r\s_t
&
~=~
(u_ru_t\otimes\id_{\breve{F}})
\begin{cases}
(\I_{\RO{0,1-r-t}}+\I_{\RO{1-t,1}})\I_{\RO{0,1-t}}\otimes\id_{\breve{F}}		&
\\
~~~+\bfam{(\I_{\RO{0,1-r-t}}+\I_{\RO{1-t,1}})\I_{\RO{1-t,1}}+\I_{\RO{1-r-t,1-t}}\I_{\RO{0,1-t}}}\otimes\breve{s}		&
\\
~~~~~~~~~~~~~~~~~~~~~~~~~~~~~~~~~~~~~~~~~~~~~~~~+\I_{\RO{1-r-t,1-t}}\I_{\RO{1-t,1}}\otimes\breve{s}^2		&~~~~~~t\le1-r,
\\
\I_{\RO{1-t,2-r-t}}\I_{\RO{0,1-t}}\otimes\id_{\breve{F}}		&
\\
~~~+\bfam{\I_{\RO{1-t,2-r-t}}\I_{\RO{1-t,1}}+(\I_{\RO{0,1-t}}+\I_{\RO{2-r-t,1}})\I_{\RO{0,1-t}}}\otimes\breve{s}		&
\\
~~~~~~~~~~~~~~~~~~~~~~~~~~~~~~~~~~~~+(\I_{\RO{0,1-t}}+\I_{\RO{2-r-t,1}})\I_{\RO{1-t,1}}\otimes\breve{s}^2		&~~~~~~t\ge1-r
\end{cases}
\\[2ex]
&
~=~
(u_{r+t}\otimes\id_{\breve{F}})
\begin{cases}
\I_{\RO{0,1-r-t}}\otimes\id_{\breve{F}}
+\bfam{\I_{\RO{1-t,1}}+\I_{\RO{1-r-t,1-t}}}\otimes\breve{s}
+0\otimes\breve{s}^2		&~~~~~~~~~~~~~~~t\le1-r,
\\
0\otimes\id_{\breve{F}}
+\bfam{\I_{\RO{1-t,2-r-t}}+\I_{\RO{0,1-t}}}\otimes\breve{s}
+\I_{\RO{2-r-t,1}}\otimes\breve{s}^2		&~~~~~~~~~~~~~~~t\ge1-r
\end{cases}
\\[2ex]
&
~=~
(u_{r+t}\otimes\id_{\breve{F}})
\begin{cases}
\I_{\RO{0,1-r-t}}\otimes\id_{\breve{F}}+\I_{\RO{1-r-t,1}}\otimes\breve{s}		&~~~~~~~~~~~~~~~~~~~~~~~~~~~~~~~~~~~~~~~~~~~~~~~~t\le1-r,
\\
\I_{\RO{0,2-r-t}}\otimes\breve{s}+\I_{\RO{2-r-t,1}}\otimes\breve{s}^2		&~~~~~~~~~~~~~~~~~~~~~~~~~~~~~~~~~~~~~~~~~~~~~~~~t\ge1-r
\end{cases}
\\
&
~=~
\s_{r+t}
}\ealn
in either case.

Part \ref{2} is shown best by decomposing $L^2(\R_+)\otimes F$ into a direct sum over $L^2\RO{n-1,n}\otimes F$ $(n\in\N)$ in which, then, each $L^2\RO{n-1,n}\otimes F$ is identified with $L^2\RO{0,1}\otimes F$. As tracking what the standard right shift on $L^2(\R_+)\otimes F$ does in this identification, is quite obvious, we do not give details.\qed

\bex\label{nondecex}
The first part of the lemma promises to find an example of a strongly continuous semigroup of adjointable isometries $\s_t$ on a Hilbert module $E$ with $E_p^\perp=\zero$ but $E_p\ne E$, provided we find an adjointable isometry $\breve{\s}$ on a Hilbert module $\breve{E}$ such that $\breve{E}_p:=\ol{\bigcup_n(\breve{\s}^n\breve{E})^\perp}\ne\breve{E}$ but $\breve{E}_p^\perp\ne\zero$. (Indeed, defining $E$ and $\s_t$ as in the lemma, we conclude from $\s_n=\id\otimes\breve{\s}^n$ and from $\s_t\s_t^*$ increasing, that $E_p=L^2\RO{0,1}\otimes\breve{E}^p\ne E$. The statement about the complement is shown precisely as for $\cI$ and $\cB$ in the beginning of Example \ref{nonadex}.)

Let $\breve{v}$ be a proper isometry on a Hilbert space $H\ne\zero$. Define $\breve{E}:=C_b(\N,H)$, the bounded \nbd{H}valued functions on $\N$. It is routine to show that the inner product $\AB{f,g}$ defined by setting $\AB{f,g}(k):=\AB{f(k),g(k)}$ turns $\breve{E}$ into a Hilbert module over $\cB:=C_b(\N)$. Then the operator $\breve{\s}$ on $\breve{E}$ defined by pointwise action of $\breve{v}$ on a function $f\in\breve{E}$, $\SB{\breve{\s}f}(k):=\breve{v}f(k)$, is an isometry and pointwise action of $\breve{v}^*$ is an adjoint. Since $\breve{v}$ is proper, we may choose an orthonormal family $\CB{e_n}_{n\in\N}$ in $H$ such that $e_n\in(\breve{v}^{n-1}\breve{v}^{*^{\scriptstyle n-1}}-\breve{v}^n\breve{v}^{*^{\scriptstyle n}})H$. Define $f\in\breve{E}$ by $f(k)=e_k$. Then, clearly, $\lim_{n\to\infty}\breve{v}^n\breve{v}^{*^{\scriptstyle n}}f$ does not exists, so $f\notin\breve{E}_p$. On the other hand, if $\breve{v}$ is pure (that is $\breve{v}^{*^{\scriptstyle n}}$ converges strongly to $0$), then, clearly, $\breve{E}_p^\perp=\zero$.
\eex

\bob
Note that in the preceding example with pure $\breve{v}$, it is easy to see that $E_p$ may be identified with $H\otimes\cB$. Since every infinite-dimensional Hilbert space admits a pure isometry, this also shows that in this case $C_b(H)\supsetneq H\otimes C_b(\N)$.

One may ask, why in the example, not taking immediately a pure strongly continuous semigroup of isometries $\breve{v}_t$ on $H$ (necessarily unitarily equivalent to the standard right shift on some $L^2(\R_+,K)$)? In fact, there is no problem to define, then, on the same $\breve{E}=C_b(\N,H)$ a semigroup of adjointable isometries $\breve{\s}_t$ by pointwise action of $\breve{v}_t$ on functions on $\breve{E}$. However, the action of $\breve{\s}$ on the function $g\in\breve{E}$ defined by
\beqn{
g(k)
~:=~
y\frac{\I_{\frac{1}{k+1},\frac{1}{k}}}{\sqrt{\frac{1}{k}-\frac{1}{k+1}}}
~=~
y\sqrt{k(k+1)}\I_{\frac{1}{k+1},\frac{1}{k}}
}\eeqn
for some non-zero vector $y\in K$, shows that the semigroup $\breve{\s}_t$ is not strongly continuous. It is strongly continuous when restricted to the submodule $\breve{E}_p=H\otimes C_b(\N)\notni g$. The same $g$ shows that $\breve{\s}_t$ is not even \nbd{\cB}weakly continuous, that is, $t\mapsto\AB{g,\breve{\s}_tg}$ is not continuous (in any standard topology of $\cB$).
\eob

\bob\label{notscob}
Note that every semigroup of (not necessarily adjointable) isometries $\s_t$ on $E$ may be `dilated' to a semigroup $\wt{\s}_t$ of unitaries on $\wt{E}$. (`Dilating' means that $\wt{E}\supset E$ and that $\wt{\s}_t$ (co)restricts to $\s_t$ on $E$.) Indeed, if we put $E_t:=E$ and define the maps $\beta_{t,s}\colon E_s\rightarrow E_t$ for all $t\ge s$ to be $\s_{t-s}$, then the $E_t$ and $\beta_{t,s}$ form an inductive system. Let $\wt{E}$ denote the inductive limit and denote by $i_t$ the canonical embeddings of $E_t$ into $\wt{E}$. It is not difficult to show that $\wt{\s}_t\colon k_{s+t}x\mapsto k_sx$ determines a unitary for each $t$ and that these unitaries form a semigroup $\wt{\s}$. (See the appendix of Bhat and Skeide \cite{BhSk00} for inductive limits of Hilbert modules.) Of course, $E\cong i_0E\subset\wt{E}$ and restriction to $E$ gives back $\s$. (The inductive limit for the semigroup of non necessarily adjointable isometries in Example \ref{nonadex} is $L^2(\R,\cB)$ with canonical embeddings $i_t$ sending $L^2(\R_-,\cI)\oplus L^2(\R_+,\cB)$ onto the subset $L^2(\LO{-\infty,-t},\cI)\oplus L^2(\RO{-t,\infty},\cB)$ of $L^2(\R,\cB)$.)

One may show: $\wt{\s}$ is strongly continuous if and only if $\s$ is; $E$ is complemented in $\wt{E}$ if and only if the $\beta_{t,s}$ are adjointable, that is, if the $\s_t$ are adjointable; if $\s$ is even pure, then the projections $p_t$ onto $(\wt{\s}_tE)^\perp$ converge strongly to $\id_{\wt{E}}$ and $0$ for $t\to\infty$ and $t\to-\infty$, respectively. So, if $\s$ is a strongly continuous semigroup of adjointable isometries, then the family $p_\lambda$ is a spectral measure that is strongly continuous. So, it follows that the integrals
\beqn{
\wt{\st}_t
~:=~
\int e^{it\lambda}\,dp_\lambda
}\eeqn
exist strongly and define a strongly continuous group of unitaries on $\wt{E}$. Extending $\wt{\s}$ to negative times, the two strongly continuous unitary groups $\wt{\s}$ and $\wt{\st}$ fulfill the \hl{Weyl commutation relations}
\beqn{
\wt{\s}_s\wt{\st}_t
~=~
e^{ist}\wt{\st}_t\,\wt{\s}_s.
}\eeqn

Now, if we apply our theorem to the pure strongly continuous semigroup $\s$, so that $\s$ is given as a standard right shift on some $L^2(\R_+,F)$, then $\wt{E}=L^2(\R,F)$ with $i_tE=L^2(\RO{-t,\infty},F)$, and
\beqn{
\SB{\wt{\st}_t(f)}(x)
~=~
e^{itx}f(x).
}\eeqn
The Stone-von Neumann theorem asserts that every pair of strongly continuous semigroups of unitaries on a Hilbert space $H$ that fulfills the Weyl relations, is unitarily equivalent to the pair $\wt{\s}$ and $\wt{\st}$ for some Hilbert space $F$.

This can be interpreted in two directions. On the one hand, whenever the Stone-von Neumann theorem holds, we can use this to prove Theorem \ref{mainthm}. (Simply construct $\wt{E}$, $\wt{\s}$, and $\wt{\st}$, apply the Stone-von Neumann theorem, and remember how $E$ and $\s$ sit inside $L^2(\R,F)$.) On the other hand, if, starting with strongly continuous unitary groups $\wt{\s}$ and $\wt{\st}$ on $\wt{E}$ that fulfill the Weyl relations, we would succeed to find a submodule $E$ turning the $p_\lambda$ (defined as above) into a spectral measure giving $\wt{\st}_t$ back as $e^{it\lambda}$, then Theorem \ref{mainthm} would allow prove the Stone-von Neumann theorem for that pair.

For Hilbert spaces the latter can be done, writing $\wt{\st}_t$ as $\int_\R e^{it\lambda}\,dp_\lambda$ by means of Stone's theorem, and showing that $\wt{\s}_s\,dp_\lambda=dp_{\lambda+s}\,\wt{\s}_s$. (Then $E:=p_0$ does the job.) For Hilbert modules there is no spectral theorem, consequently, there is no Stone theorem for unitary groups. We, thus, do not know how to restrict one of the unitary groups fulfilling the Weyl relations to a pure semigroup of isometries to which our Theorem could be applied. Consequently, we do not get a Stone-von Neumann theorem for unitary groups on Hilbert modules.
\eob

\bob
Note that for von Neumann (or $W^*$-modules), Theorem \ref{SzNFthm} generalizes verbatim, if we understand that ``strongly continuous'' is now referring to the strong operator topology of the von Neumann module $E\subset\sB(G,H)$, where the von Neumann algebra $\cB$ is acting (nondegenerately) on a Hilbert space $G$, and where the Hilbert space $H=E\odot G$ is the internal tensor product over $\cB$ of $E$ with the Hilbert \nbd{\cB}\nbd{\C}module $G$; see Skeide \cite{Ske00b,Ske01}. (The fact that the semigroup in Theorem \ref{mainthm} is a $C_0$-semigroup, which is much stronger a continuity condition, plays a crucial role in the proof in Section \ref{proofSEC}. Like all bounded right linear maps on a von Neumann module, the $\s_t$ are adjointable, automatically. Also, $E_u^\perp=\ol{E_p}^s$ (the strong closure of $E_p$) and $E=E_u\oplus E_u^\perp$. In the case $E=E_u^\perp$ (equivalently, $E_u=\zero$), necessarily $\s_t\s_t^*\uparrow\id_E$ in the strong topology of $\sB^a(E)\subset\sB(H)$. We may apply Theorem \ref{SzNFthm} to that semigroup on $H$, and obtain a Hilbert space $K$ such that $\s_t$ is given by $v_t$ on $L^2(\R_+,K)$.

Now by Skeide \cite{Ske05c}, the subspace $E\subset\sB(G,H)$ is characterized as
$$
E
~=~
\bCB{x\in\sB(G,H)\colon \rho'(b')x=xb'~(b'\in\cB')},
$$
where $\rho'$ is a (unique) normal unital representation of $\cB'$ on $H$, and $\sB^a(E)=\rho'(\cB')'$. Since $H=L^2(\R_+,K)$ and $\s_t\in\sB^a(E)$, there is a representation of $\cB'$ on $L^2(\R_+,K)=L^2(\R_+)\otimes K$ such that its range commutes with all $v_t$. It is standard to verify, that this representation is necessarily of the form $\id_{L^2(\R_+)}\otimes\sigma'$, where $\sigma'$ is a (normal unital) representation of $\cB'$ on $K$. By \cite{Ske03c}, the space
$$
F
~:=~
\bCB{y\in\sB(G,K)\colon \sigma'(b')x=xb'~(b'\in\cB')}
$$
is a von Neumann $\cB$-module. It is routine to verify that the strong closure of $L^2(\R_+)\otimes F$ in $\sB(G,L^2(\R_+)\otimes K)$ (the element $f\otimes y$ acting as $g\mapsto f\otimes yg$) is $E$.

Of course, in the very same way we do get a Stone-von Neumann theorem for strongly continuous semigroups of unitaries on a von Neumann module.
\eob

\section{Proof of Theorem \ref{mainthm}} \label{proofSEC}

Let $\cB$ be a  $C^*$-algebra and let $E$ be a Hilbert
$C^*$-module over $\cB$. Let $\sB^a(E)$ denote the $C^*$-algebra of all
adjointable (hence bounded) operators on $E$.  Let $\{ \s_t\}_{t\geq
0}$ be a strongly continuous one parameter semigroup of adjointable
isometries on $E$, that is,

Note that $\s_t\s_t^*$ is a decreasing family of projections with
$\s_0\s_0^*=\id$. So for $0\leq a <b< \infty$ also $p_{a,b}:=
\s_a\s_a^*-\s_b\s_b^*$ form a spectral measure on the reals half-line. We also put $p_{c,d}:=0$ if
$c>d$, and $p_{a, \infty}:=\s_a\s_a^*$. We take $E_{a,b}:= p_{a,b}E$.

Recall the standard notation $c\vee d:= \max \{ c , d\},$ and $c\wedge d := \min  \{ c , d\}.$

\begin{prop} \label{pabprop}
Let $0\leq t<\infty,$ $0\leq a<b\le\infty$ and $0\leq c<d\le\infty$.

(i) $p_{a,b}p_{c,d}= p_{a\vee c, b\wedge d}$.

(ii) $\s_tp_{a,b}= p_{a+t , b+t}\s_t$; $p_{a,b}\s_t=\s_tp_{(a-t)\vee 0,
(b-t)\vee 0}$

(iii) $ \s_t^*p_{a,b}=p_{(a-t)\vee 0,(b-t)\vee 0 }\s_t^*$;
$p_{a,b}\s_t^*=\s_t^*p_{a+t, b+t}$.
\end{prop}

\proof (i) is a standard computation for spectral measures.

The first formula in (ii) follows from 
$\s_t\s_a\s_a^*= \s_{a+t}\s_a^*= \s_{a+t}\s_{a+t}^*\s_t$. The second formula  in (ii) follows by taking also into account that for $t\ge a$ we get $\s_a\s_a^*\s_t=\s_a\s_{t-a}=\s_t=\s_0\s_0^*\s_t$.

The formulae in (iii) are adjoints of (ii).\qed

\lf\noindent
The $E_{a,b}$ are our candidates for $L^2(\RO{a,b},F)$. A typical behaviour is:

\bcor \label{abshiftcor}
From $\s_tp_{a,b}= p_{a+t , b+t}\s_t$ and $p_{a,b}\s_t^*=\s_t^*p_{a+t, b+t}$ we infer that $\s_t$ (co)restricts to a unitary $E_{a,b}\rightarrow E_{a+t,b+t}$ with inverse $\s_t^*$ restricted to $E_{a+t,b+t}$.
\ecor

\lf
We are now going to define the unitary groups $u^{a.b}$ that simulate the shift modulo $b-a$ on $L^2(\RO{a,b},F)$. For $0\leq a<b<\infty $ and $0\leq t < (b-a)$, define $u_t^{a,b}$ by
$$u_t^{a,b}= \s_tp_{a, b-t}+\s_{b-a-t}^*p_{b-t, b}.$$
Making use of Proposition \ref{pabprop}, we also have:
\beq{\label{uabdefi}
u_t^{a,b}= p_{a+t, b}\s_t+p_{a, a+t}\s_{b-a-t}^*.
}\eeq
We extend the definition of $u^{a,b}_t$  periodically to all $t\in \mathbb{R}$
by setting
$$u^{a,b}_t = u^{a,b}_{t-n(b-a)},$$
for $t\in [n(b-a), (n+1)(b-a))$, with  $n\in \mathbb{Z}.$

\begin{prop} (Unitarity and Group property)
Each $u_t^{a,b}$ is a unitary on $E_{a,b}$ and $u_r^{a,b}u_t^{a,b}= u_{r+t}^{a,b}$
for all $ r, t\in \mathbb{R} $.

\end{prop}


\proof 
Since $u_t$ is periodic with period $b-a$, it is sufficient if we do computations for times in $\RO{0,b-a}$.

By Corollary \ref{abshiftcor}, $u^{a,b}_t$, as the direct sum of two unitaries from $E_{a,b}=E_{a, b-t}\oplus E_{b-t, b}$ onto $E_{a+t, b}\oplus E_{a, a+t}=E_{a,b}$, is unitary.

For the semigroup property, let $0\leq r, t<b-a$. Then
$$\begin{array}{rcl}
u_r^{a,b}u_t^{a,b} &=& (\s_rp_{a, b-r}+\s_{b-a-r}^*p_{b-r, b})
(p_{a+t, b}\s_t+p_{a, a+t}\s_{b-a-t}^*)\\
&=& \s_rp_{a\vee (a+t), (b-r)\wedge b}\s_t+ \s_rp_{a, (b-r)\wedge
(a+t)}\s_{b-a-t}^*\\
&&~~~~~~~~~~~~~~~~~~~~~~~~~~ + \s_{b-a-r}^*p_{(b-r)\vee (a+t), b}\s_t+
\s_{b-a-r}^*p_{b-r\vee a, b\wedge (a+t)}\s_{b-a-t}^*\\
&=& \s_rp_{a+t, b-r}\s_t+ \s_rp_{a, (b-r)\wedge
(a+t)}\s_{b-a-t}^*+ \s_{b-a-r}^*p_{(b-r)\vee (a+t), b}\s_t+
\s_{b-a-r}^*p_{b-r, a+t}\s_{b-a-t}^*.\\
\end{array}$$
We have to distinguish two cases. Firstly, $a+t<b-r$, that is, $r+t<b-a$. So,
$$\begin{array}{rcl}
u_r^{a,b}u_t^{a,b} &=&\s_rp_{a+t, b-r}\s_t+ \s_rp_{a, a+t}\s_{b-a-t}^*+ \s_{b-a-r}^*p_{b-r,
b}\s_t+ 0\\
&=& \s_{r+t}p_{a, b-r-t}+ p_{a+r, a+r+t}\s_r\s_{b-a-t}^*+
\s_{b-a-r}^*\s_tp_{b-r-t, b-t}\\
&=& \s_{r+t}p_{a, b-r-t}+ p_{a+r, a+r+t}p_{r, \infty}\s_{b-a-r-t}^*+
\s_{b-a-r-t}^*p_{b-r-t, b-t}\\
&=& \s_{r+t}p_{a, b-r-t}+ p_{a+r, a+r+t}\s_{b-a-r-t}^*+
\s_{b-a-r-t}^*p_{b-r-t, b-t}\\
&=& \s_{r+t}p_{a, b-r-t}+ \s_{b-a-r-t}^*p_{b-t, b}+
\s_{b-a-r-t}^*p_{b-r-t, b-t} \\
&=& \s_{r+t}p_{a, b-r-t}+ \s_{b-a-r-t}^*p_{b-r-t, b}\\
&=& u^{a,b}_{r+t}.\\
\end{array}$$

Secondly, $a+t\ge b-r$, that is, $r+t\ge b-a$. So,
$$\begin{array}{rcl}
 u_r^{a,b}u_t^{a,b} &=&  0+ \s_rp_{a, b-r}s^*_{b-a-t} + s^*_{b-a-r}p_{a+t,
b}\s_t+s^*_{b-a-r}p_{b-r, a+t}s^*_{b-a-t}\\
&=& \s_{r+t-(b-a)}p_{b-t, 2b-a-r-t}+\s_{t-(b-a-r)}p_{a,
b-t}+s^*_{2(b-a)-r-t}p_{2b-a-r-t, b}\\
&=& \s_{r+t-(b-a)}p_{a, 2b-a-r-t}+ s^*_{2(b-a)-r-t}p_{2b-a-r-t, b}\\
&=& u^{a,b}_{r+t-(b-a)}~=~u^{a,b}_{r+t}.\\
\end{array}$$
Now, a semigroup homomorphism between groups is a group homomorphisms. So, the $u^{a,b}_t$, effectively, form a unitary one-parameter group.\qed

\begin{prop} (Continuity) $ t\mapsto u^{a,b}_t$ is
strongly continuous.

\end{prop}

\proof Observe that for $x\in E$ and $t\geq 0$, $\| \s_t^*x-x\|= \|
\s_t^*x-\s_t^*\s_tx\| \leq \| \s_t^*\| (\| x-\s_tx\| )= \| x-\s_tx\|.$
Hence $t\mapsto \s_t^*$ is strongly continuous. If $t\mapsto A_t$ and $t\mapsto B_t$ are strongly continuous, then $t\mapsto A_tB_t$ is
also strongly continuous. (Indeed, by the \it{principle of uniform boundedness}, $A_t$ is locally uniformly bounded. By $\|A_tB_tx-A_rB_rx\|\leq \|A_t\| \,\|
(B_t-B_r)x\| + \|(A_t-A_r)B_rx\|$, we see that for fixed $r\in\R$ and $x\in E$, $A_tB_tx$ converges to $A_rB_rx$ for $t\to r$.) Therefore, maps like $t\mapsto
\s_{t+\alpha}\s_{t+\beta}^*$ and $t\mapsto
\s_{t+\alpha}^*\s_{t+\beta}$ are strongly continuous. Consequently, $t\mapsto u^{a,b}_t$
is strongly continuous on $[0, (b-a))$. Further,
$$\begin{array}{rcl}
\displaystyle\lim_{t\uparrow (b-a)} u^{a,b}_t &=& \displaystyle\lim_{t\uparrow
(b-a)}[\s_tp_{a,b-t}+s^*_{b-a-t}p_{b-t, b}]\\
&=& \displaystyle\lim_{t\uparrow
(b-a)}[\s_t(\s_a\s_a^*-\s_{b-t}\s_{b-t}^*)+s^*_{b-a-t}(\s_{b-t}\s_{b-t}^*-\s_b\s_b^*)]\\
&=&0+\id(\s_a\s_a^*-\s_b\s_b^*) =p_{a,b}=u_0^{a,b}.
\end{array} $$
Therefore, also the periodic extension to the real line is strongly continuous.\qed

\lf
The functions of the form $\I_{\RO{a,b}}y$ in $L^2(\RO{a,b},F)$ are precisely those which are invariant under the unitary shift modulo $b-a$. Taking the mean should give as a projection onto that invariant subspace. Therefore, for $x\in E $, and $0\leq  a< b < \infty $, define
$$q_{a,b}x= \frac{1}{b-a}\int _0^{b-a}u_t^{a,b}x\,dt$$
in the sense of Riemann integral over the continuous function $u^{a,b}_tx$.

\blem
$u^{a,b}_rq_{a,b}=q_{a,b}$.
\elem

\proof
It is enough to show the statement for $0\leq r<(b-a)$. As $u_t^{a,b} = u_{t-(b-a)}^{a,b}$ for $b-a\leq t<2(b-a)$, we get
$$\begin{array}{rcl}
u_r^{a,b}q_{a,b}x &=&\displaystyle \frac {1}{b-a}\int _0^{b-a}
u_r^{a,b}u_t^{a,b}x\,dt\\
&=&\displaystyle \frac {1}{b-a}\int _r^{(b-a)+r}u_t^{a,b}x\,dt\\
&=&\displaystyle \frac {1}{b-a}\int _r^{b-a}u_t^{a,b}x\,dt+\frac {1}{b-a}\int
_0^ru_t^{a,b}x\,dt \\
&=& q_{a,b}x.\qedsymbol
\end{array}$$
\noqed

\bcor\label{subpcor}
$q_{a,b}$ is a subprojection of $p_{a,b}$ on
$E$.
\ecor

\proof
As $$u_t^{a,b}= p_{a+t, b}\s_t+p_{a, a+t}\s_{b-a-t}^*,$$ it is clear
that the range of $q_{a,b}$ is contained in the range of $p_{a,b}$.

To show that a linear map $q$ is a projection, it is sufficient to check $\AB{qx,qy}=\AB{x,qy}$. For $q_{a,b}$ that formula follows from $\AB{qx,qy}=\frac {1}{b-a}\int _0^{b-a}\AB{x,(u^{a,b}_t)^*q_{a,b}y}$ and the lemma.\qed

\lf
The following simple corollary of the lemma shows that pieces of $q_{a,b}x$ in subintervals behave nicely with respect to the shift.

\bcor\label{shift}
For $0<r<b-a, $,
$$ \s_rp_{a , b-r}q_{a,b}= p_{a+r , b}\s_rq_{a,b}= p_{a+r , b}u^{a,b}_rq_{a,b}=p_{a+r, b}q_{a,b}.$$
\ecor

\noindent
We now come to the analytical heart of our proof of Theorem \ref{mainthm}. The following lemma will guarantee that we may approximate everything by linear combinations of elements $q_{a,b}x$ ($x\in E$, $0\le a<b<\infty$).

\blem\label{limit}
For $x\in E_{0,1},$
$$
\lim _{n\to \infty }\sum _{k=1}^n
q_{\frac {k-1}{n}, \frac{k}{n}}x =x.
$$
\elem

\proof
Note that
$$
\sum _{k=1}^n
q_{\frac {k-1}{n}, \frac{k}{n}}x =
 n\int _0^{\frac{1}{n} }\sum _{k=1}^n
u_t^{\frac {k-1}{n}, \frac{k}{n}}x\,dt.
$$
For $\epsilon >0$, choose $N$ such that  for all $t\in\bSB{0,
\frac{1}{N}},$
$$\| \s_tx-x\| <\frac{\epsilon }{2}, ~~\mbox{and} ~~ \| \s_t^*x-x\|
<\frac{\epsilon }{2}.$$ Consider $n\geq N$. We get
\baln{
\Bnorm{ \Bfam{\sum _{k=1}^n u_t^{\frac {k-1}{n}, \frac{k}{n}} x}-x}
& =
\Bnorm{\Bfam{\sum _{k=1}^np_{{\frac {k-1}{n}+t}, \frac{k}{n}}\s_tx+\sum
 _{k=1}^np_{{\frac {k-1}{n}}, \frac{k-1}{n}+t}s^*_{\frac
 {1}{n}-t}x}
 -
\Bfam{\sum _{k=1}^np_{{\frac {k-1}{n}+t}, \frac{k}{n}}x+\sum
 _{k=1}^np_{{\frac {k-1}{n}}, \frac{k-1}{n}+t}x} }
\\
 &\leq   
\Bnorm{ \sum _{k=1}^np_{{\frac {k-1}{n}+t},
 \frac{k}{n}}(\s_tx-x)}
 + 
\Bnorm{ \sum _{k=1}^n p_{{\frac {k-1}{n}},
 \frac{k-1}{n}+t} (s^*_{\frac {1}{n}-t}x-x)}
 \\
 &\leq
\Bnorm{ \sum _{k=1}^n p_{{\frac {k-1}{n}+t}, \frac{k}{n}}}\, \| \s_tx-x \|
 + 
\Bnorm{ \sum _{k=1}^n p_{{\frac {k-1}{n}},\frac{k-1}{n}+t}}\, \| (\s_{\frac{1}{n}}^*x-x)\| 
\\
  &\leq  
\frac{\epsilon }{2}+\frac {\epsilon} {2} =\epsilon .
 }\ealn
Hence,
$\bnorm{ n\int _0^{\frac{1}{n}} \sum _{k=1}^n u_t^{\frac {k-1}{n}, \frac{k}{n}}x\,dt -x}
 =
\bnorm{ n \int _0^{\frac {1}{n}} \bfam{\sum_{k=1}^n u_t^{\frac {k-1}{n}, \frac {k}{n}}-\id}x \,dt}
\leq 
 n\int _0^{\frac{1}{n}}\epsilon \,dt
 =
 \epsilon$.\qed

\lf
Form now on we prepare for proving in Proposition \ref{01isoprop} that $E_{0,1}\cong L^2(\RO{0,1},q_{0,1}E)$, making $q_{0,1}E$ our `hot' candidate for being the $F$ we seek, and for proving in Corollary \ref{intercor} that $u^{0,1}$, indeed, transfroms under this isomorphisms into the shift modulo $1$.

\begin{prop}\label{ratio}

For $0\leq a<c<d<b<\infty$,
$$p_{c,d}q_{a,b}p_{c,d}= \frac{d-c}{b-a}q_{c,d}.$$
\end{prop}
\proof
For $0\leq t\leq (b-a),$
$$\begin{array}{rcl}
p_{c,d}(p_{a+t,b}\s_t+\s_t^*p_{a+t,b})p_{c,d} &=&\displaystyle p_{c\vee a+t,
d}p_{c+t, d+t}\s_t+ \s_t^*p_{c+t, d+t}p_{c\vee a+t, d} \\
&=&\displaystyle p_{c+t, d}\s_t+\s_t^*p_{c+t, d},\end{array}$$
 which is non-zero only if $0< t<(d-c).$ Hence, by an the application of change of variable
$$\begin{array}{rcl}
q_{a,b}x&=&  \displaystyle\frac{1}{b-a}\int _0^{b-a} (p_{a+t, b}\s_t+p_{a,
a+t}\s_{b-a-t}^*)x\,dt\\
&= &\displaystyle
 \frac{1}{b-a}\Bfam{\int _0^{b-a} p_{a+t, b}\s_tx\,dt+ \int_0^{b-a}p_{a,b-t}\s_t^*x\,dt}  ~~~~ (z=b-a-t)\\
 &=&\displaystyle \frac{1}{b-a}\int _0^{b-a} (p_{a+t, b}\s_t+ \s_t^*p_{a+t,b})x\,dt,
 \end{array}$$
we get
 $$\begin{array}{rcl}
 p_{c,d}q_{a,b}p_{c,d}x &=&\displaystyle \frac{1}{b-a}\int
 _0^{d-c}(p_{c+t,d}\s_t+\s_t^*p_{c+t,d})x\,dt\\
 &=&\displaystyle \frac{d-c}{b-a}q_{c,d}x.\qedsymbol\end{array}$$\noqed

\begin{prop}
$\s_rq_{a,b}= q_{a+r, b+r}\s_r$, that is, $\s_r$ (co)restricts to a unitary $q_{a,b}E\rightarrow q_{a+r,
b+r}E$.
\end{prop}

\proof
$$\begin{array}{rcl}
\s_rq_{a,b}x &= &\displaystyle \frac{1}{b-a}\int _0^{b-a}\s_r[\s_tp_{a, b-t}+p_{a,a+t} \s_{b-a-t}^*]x \,dt \\
&=&\displaystyle \frac{1}{b-a}\int _0^{b-a} [\s_tp_{a+r, b+r-t}\s_rx+p_{a+r,a+r+t}\s_rs^*_{b-a-t}x]\,dt\\
&=&\displaystyle \frac{1}{b-a} \int _0^{b-a}[\s_tp_{a+r, b+r-t}\s_rx+p_{a+r, a+r+t}p_{r, \infty }s^*_{b-a-t}\s_rx]\,dt\\
&=&\displaystyle \frac{1}{b-a} \int _0^{b-a} [\s_tp_{a+r, b+r-t}+p_{a+r, a+r+t} s^*_{(b+r)-(a+r)-t}]\s_rx\,dt\\
&=& q_{a+r, b+r}\s_rx.\qedsymbol
\end{array}$$ \noqed

\begin{prop}
Suppose $z\in q_{0,1}(E)$. Then for $0<r,t$, with $r+t<1$,
$$p_{0, r+t}z=p_{0,r}z+\s_rp_{0,t}z.$$

\end{prop}

\proof
We have $ p_{0, r+t}z=p_{0,r}z+p_{r,r+t}z= p_{0,r}z+p_{r,
r+t}p_{r,1}z.$ Now   Proposition\ref{ratio}, with $a=0, b=1$ yields
$p_{r,1}z=p_{r,1}\s_rz.$  Hence, $p_{0, r+t}{z}= p_{0,r}z+p_{r,
r+t}p_{r,1}\s_rz= p_{0,r}z+p_{r, r+t}\s_rz=p_{0,r}z+\s_rp_{0,t}z.$ \qed


\begin{prop}
Suppose $z, w\in q_{0,1}(E)$. Then for $0<c_1<d_1<1$ and
$0<c_2<d_2<1$, $\langle p_{c_1, d_1}z, p_{c_2, d_2}w\rangle= \mu
([c_1, d_1]\bigcap [c_2 , d_2])\langle z, w\rangle $, where $\mu $
denotes the Lebesgue measure.
\end{prop}
\proof
For $0<t<1$, set $f(t)= \langle z , p_{0, t}w\rangle .$ Then from
the previous Proposition, for $0<r+t<1$,
$$\begin{array}{rcl}
f(r+t) &=& \langle p_{0, r+t}z , p_{0, r+t}w\rangle \\
&=& \langle  p_{0, r}z ,p_{0,r}w\rangle +\langle \s_rp_{0,t}z,
\s_rp_{0,t}w\rangle \\
&=& \langle p_{0,r}z, p_{0,r}w\rangle + \langle p_{0, t}z,
p_{0,t}w\rangle \\
&=& f(r)+f(t).
\end{array} $$  Then by strong continuity of $t\mapsto p_{0
,t}$, it follows that $f(t)= t\langle z , w\rangle $ for all $t$.
Hence the result. \qed

\lf
Recall that  $L^2([0 , 1], F):=L^2[0,1]\otimes F$ is the external tensor product. The subset of functions of the form $\I _{[c , d)}z$ ($0\leq c<d\leq 1$) is total.

\begin{prop} \label{01isoprop}
Take $F=q_{0,1}(E).$ Define $M: L^2([0 , 1], F)\to E_{0, 1}$ by
$$M(\I _{[c , d)}z)= p_{c,d}z,$$
for $z\in F$ and $0\leq c<d\leq 1$ ($\I $ denotes the indicator
function), and extending linearly. Then $M$ extends to a unitary
map.
\end{prop}

\proof
The isometry property of $M$ has been proved in the previous
proposition. Now for $x\in E_{0,1},$ and $1\leq k\leq n$, by
Proposition \ref{ratio}, $ q_{\frac{(k-1)}{n}, \frac{k}{n}}x$ is in
the range of $p_{\frac{(k-1)}{n}, \frac{k}{n}}q_{0,1}$.
 Then from Proposition \ref{limit},
$  \lim _{n\to \infty}\sum _{k=1}^n  q_{\frac{(k-1)}{n},
\frac{k}{n}}x   =\lim _{n\to \infty } n\int _0^{\frac{1}{n} }\sum
_{k=1}^n u_t^{\frac {k-1}{n}, \frac{k}{n}}x\,dt =x.$ This shows that
the range of $M$ is whole of $E_{0,1}.$ \qed

\lf
Let $\pi_t$ denote the periodic shift on $L^2([0 , 1], F)$

\bcor\label{intercor}
$M^*u^{0,1}_t=\pi_tM^*$.
\ecor

\proof
It suffices to check this on $p_{c,d}z$ with either $0\le c<d<1-t$ or $1-t\le c<d<1$, as every other choice is a sum of such. For the stated cases, the statements follows directly from the definitions of $u^{0,1}_t$, $\pi_t$, and $M$.\qed

\brem
Of course, the proposition and its corollary remain true, replacing the interval $\RO{0,1}$ with any other interval nonempty interval $\RO{a,b}\subset\R_+$.
\erem

\proof[Proof of Theorem \ref{mainthm}.]
By pureness,
\beqn{
E=\bigoplus _{k=0}^{\infty }E_{k,k+1}=\bigoplus _{k=0}^{\infty }\s_kE_{0,1}.
}\eeqn
Since $\s_kE_{0,1}\cong L^2\RO{0,1}\otimes F$, and since $\s$ acts under these identifications as it should, we are ready for an application of Lemma \ref{interlem}\eqref{2}.\qed

\lf\noindent
\bf{Acknowledgments.~}
We started thinking about Theorem \ref{mainthm} during a {Research in Pairs} project at {Mathematisches Forschungsinstitut Oberwolfach} in 2007; hospitality at MFO during this RiP is gratefully acknowledged. MS is grateful to BVRB and {The Indian Statistical Institute, Bangalore Center} for kind hospitality during several stays, and for travel support from the Dipartimento E.G.S.I. of University of Molise.


\newcommand{\Swap}[2]{#2#1}\newcommand{\Sort}[1]{}
\providecommand{\bysame}{\leavevmode\hbox to3em{\hrulefill}\thinspace}
\providecommand{\MR}{\relax\ifhmode\unskip\space\fi MR }
\providecommand{\MRhref}[2]{%
  \href{http://www.ams.org/mathscinet-getitem?mr=#1}{#2}
}
\providecommand{\href}[2]{#2}

\end{document}